\documentstyle[12pt,amsmath,amssymb]{article}

\newtheorem{lem}{Lemma}[section]
\newtheorem{th}{Theorem}[section]

\newtheorem{rem}{Remark}[section]

\newtheorem{cor}{Corollary}[section]

\allowdisplaybreaks

\makeatletter\@addtoreset{equation}{section}\makeatother

\newcommand{\N}{{\Bbb N}}
\newcommand{\C}{{\Bbb C}}
\newcommand{\Z}{{\Bbb Z}}
\newcommand{\R}{{\Bbb R}}

\newcommand\lmod{\left\vert}
\newcommand\rmod{\right\vert}
\newcommand{\md}[1]{{\lmod #1 \rmod}}

\newcommand{\sk}[1]{{\left( #1 \right)}}
\newcommand{\skf}[1]{{\left\{ #1 \right\}}}

\newcommand{\vol}{\operatorname{vol}}
\newcommand{\eps}{\varepsilon}
\newcommand{\la}{\langle}
\newcommand{\ra}{\rangle}
\newcommand{\rom}[1]{{\rm #1}}
\newcommand{\dd}{\overset{{.}{.}}}

\begin{document}

\begin{center}{\Large \bf
 Non-equilibrium stochastic dynamics in continuum: The free case
}\end{center}

{\large Yuri Kondratiev}\\
 Fakult\"at f\"ur Mathematik, Universit\"at
Bielefeld, Postfach 10 01 31, D-33501 Bielefeld, Germany;
Institute of Mathematics, Kiev, Ukraine; BiBoS, Univ.\ Bielefeld,
Germany\\ e-mail:
\texttt{kondrat@mathematik.uni-bielefeld.de}\vspace{2mm}

{\large Eugene Lytvynov}\\ Department of Mathematics,
University of Wales Swansea, Singleton Park, Swansea SA2 8PP, U.K.; BiBoS, Univ.\ Bielefeld,
Germany\\
e-mail: \texttt{e.lytvynov@swansea.ac.uk}\vspace{2mm}

{\large Michael R\"ockner}\\
 Fakult\"at f\"ur Mathematik, Universit\"at
Bielefeld, Postfach 10 01 31, D-33501 Bielefeld, Germany;
 BiBoS, Univ.\ Bielefeld,
Germany\\ e-mail: \texttt{roeckner@mathematik.uni-bielefeld.de }

{\small
\begin{center}
{\bf Abstract}
\end{center}

\noindent We study the problem of identification of a proper state-space for
the stochastic dynamics of free particles in continuum, with their possible
birth and death.    In this dynamics, the motion of  each separate particle is
described by a  fixed Markov process $M$ on a  Riemannian manifold $X$. The
main problem arising here is a possible collapse  of the system, in the sense
that, though the initial configuration of particles is locally finite, there
could exist a compact set in $X$ such that, with probability one,   infinitely
many particles  will arrive  at this set at some time $t>0$.  We assume that
$X$ has infinite volume and, for each $\alpha\ge1$, we consider the set
$\Theta_\alpha$ of all infinite configurations in $X$ for which the number of
particles in a compact set is bounded by a constant times the $\alpha$-th power
of the volume of the set. We find quite general conditions on the process $M$
which guarantee that the corresponding infinite particle process  can start at
each configuration from $\Theta_\alpha$, will never leave $\Theta_\alpha$, and
has cadlag  (or, even, continuous) sample paths in the vague topology. We
consider the following examples of applications of our results: Brownian motion
on the configuration space, free Glauber dynamics on the configuration space
(or a birth-and-death process in $X$), and free Kawasaki dynamics on the
configuration space. We also show that if $X=\R^d$, then  for a wide class of
starting distributions, the (non-equilibrium) free Glauber dynamics is a
scaling limit of (non-equilibrium) free Kawasaki dynamics.
 }
 \vspace{3mm}

  \vspace{1.5mm}

\vspace{1.5mm}

\section{Introduction}  In this paper, we study the problem of identification of a proper state-space for the stochastic
dynamics of free particles in continuum, with their possible birth and death.
In this dynamics,
the motion of  each separate particle is described by a  fixed Markov process $M$ on a  Riemannian manifold
$X$.

A classical result by J.~L.~Doob \cite{Doob} states that, if the initial distribution  of free particles is Poissonian, then it will
remain Poissonian at any moment of time $t>0$, see also \cite{Dobrushin}.
 In \cite{Surgailis1,Surgailis2}, D. Surgails studied an independent infinite
particle process as a Markov process whose generator is the second quantization of  the generator of a Markov process in $X$. After Surgailis' papers, equilibrium independent infinite particle processes have been studied by many authors, see e.g.\  \cite{ST}. However, the problem of identification of allowed initial configurations of the system was not addressed in these
papers.

Let us explain this problem, in more detail. When speaking about an infinite system of particles in continuum, we should
consider such a system as an element of the configuration space $\Gamma$ over $X$. This space is defined as the
collection of all locally finite subsets of $X$. Now, consider, for example, a system of independent Brownian particles in $\R^d$. Then, one can easily find a configuration $\gamma\in\Gamma$ such that, if $\gamma$ is the initial configuration of Brownian particles, then at some time $t>0$, with probability one, the system will collapse, in the sense that there will
be an infinite number of particles in a compact set, i.e., the system will not be a configuration anymore. Thus,
generally speaking, the configuration space $\Gamma$ appears to be too big and cannot serve as a state-space for the process.

Thus, for study of the non-equilibrium stochastic dynamics of independent Markov\-ian particles, one needs to identify a subset
$\Theta$ of $\Gamma$ such that the process starting at $\Theta$ will always remain in $\Theta$ with probability one.
(We note, however, that, for a fixed independent infinite particle process, such a set $\Theta$ is not uniquely defined.)
Next, if the underlying Markov process $M$ has cadlag paths, then it is natural to expect that the corresponding infinite particle process also has cadlag sample paths in $\Theta$ with respect to the vague topology. Again, this problem was not addressed in the above-mentioned papers.

In our previous paper \cite{KLR}, we considered the case where $M$ is a
 Brownian motion in
a complete, connected, oriented, stochastically complete Riemannian
manifold $X$ of dimension $\ge2$. We explicitly constructed a subset $\Gamma_\infty$
of the configuration space and proved that the corresponding infinite particle process
  can start at any $\gamma\in\Gamma_\infty$, will never
leave $\Gamma_\infty$, and has continuous sample paths in the
vague topology (and even in a stronger one). In the case of a
one-dimensional underlying manifold $X$, one cannot exclude
collisions of particles, so that a modification of the
construction of $\Gamma_\infty$ is necessary,  see \cite{KLR} for
details.

As we mentioned above, the aim of this paper is to consider the case of general Mar\-kovian particles.
The interest in such particles, rather than just independent Brownian motions, is, in particular, connected
with study of the  Glauber and Kawasaki dynamics on the configuration space, see \cite{BCC,KL,KLR1,BCDP}.

For simpler notations, we assume that the underlying process $M$ is symmetric, however, this condition can  be easily  omitted. Instead of trying to generalize (quite complicated) arguments of \cite{KLR}, we propose a new, simpler approach
to the construction of a state-space of the
 system (which was $\Gamma_\infty$ in \cite{KLR}),
  and to the proof that the process indeed
has the above discussed properties.

So, we fix a system of closed balls $B(r)$ of radius $r\in\N$, centered at some $x_0\in X$, and for each $\alpha\ge1$,
 define $\Theta_\alpha$ as the set of those infinite configurations for which the number of particles
in each $B(r)$ is bounded by a constant times the $\alpha$-th power of the
volume of $B(r)$. The $\Theta_\alpha$'s  form an isotone sequence of sets, and
we also define $\Theta$ as the union of all the $\Theta_\alpha$'s. Note that
the sets $\Theta_\alpha$ are ``big enough'', in the sense that already the
smallest set $\Theta_1$ is of full Poisson measure with any intensity parameter
$z>0$.

We find quite general conditions on a Markov process $M$ which guarantee that the corresponding stochastic dynamics can start at each $\gamma\in\Theta_\alpha$ (or even at each $\gamma\in\Theta$), will never
leave $\Theta_\alpha$ (respectively $\Theta$), and has cadlag  (or, even, continuous) sample paths.  We then consider the following examples of application of our results:
\begin{itemize}

 \item Brownian motion on the configuration space, i.e., the case where the underlying process $M$ is a Brownian motion on
 $X$. Compared with  \cite{KLR}, our result here also covers the case of  manifolds which are not stochastically complete;

 \item Free Glauber dynamics on the configuration space, or a birth-and-death process in $X$, compare with \cite{KL,KLR1,BCDP};

 \item Free Kawasaki dynamics on the configuration space, cf.\  \cite{KLR1}.

  \end{itemize}

We also show that, in the case where  $X=\R^d$, for a probability measure $\mu$
on $\Theta$ which is translation invariant and has integrable Ursell functions,
the (non-equilibrium) free Glauber dynamics having $\mu$ as initial
distribution may be approximated by (non-equilibrium)
 Kawasaki dynamics. Note that, in the case of equilibrium dynamics,
 such an approximation can also be shown for interacting  particles, see \cite{FKL}.

  The paper is organized as follows.   In Section \ref{lkughihgh},
  we construct a Markov semigroup of kernels on the space
  $\Theta_\alpha$ and prove that it corresponds
  to the second quantization of the \linebreak (sub-)Mar\-kovian
   semigroup of the underlying process  $M$. In Section~\ref{hgig} we
   derive conditions which guarantee that the corresponding infinite
    particle process can be realized as a Markov process on $\Theta_\alpha$
    with cadlag
  (respectively, continuous) sample paths. In Section~\ref{kugyug},
  we discuss the above mentioned examples. Finally, in Section~\ref{gggh},
we prove the result on approximation of the Glauber dynamics by the Kawasaki
dynamics.

  We would like to stress that this paper should be considered as a first step towards a construction of a non-equilibrium dynamics
   of interacting particles, compare with \cite{AKR4,MR98,KL,KLR1}, where the corresponding equilibrium processes were discussed.

\section{Markov semigroup of kernels for the stochastic dynamics}\label{lkughihgh}

Let $X$ be a complete, connected, oriented $C^\infty$ Riemannian
manifold.  Let ${\cal B}(X)$ denote the Borel $\sigma$-algebra on $X$. Let $dx$ denote the volume measure on $X$, and we
suppose that $\int_Xdx=\infty$.

The configuration space $\Gamma$ over
$X$ is defined as the set of all  infinite  subsets of $X$ which are locally
finite: $$\Gamma:=\{\gamma\subset X\mid |\gamma|=\infty,\  |\gamma_\Lambda|<\infty\text{ for each compact
}\Lambda\subset X\}.$$ Here, $|\cdot|$ denotes the cardinality of
a set and $\gamma_\Lambda:= \gamma\cap\Lambda$. One can identify
any $\gamma\in\Gamma$ with the positive Radon measure
$\sum_{x\in\gamma}\eps_x\in{\cal M}(X)$. Here, $\varepsilon_x$ denotes the Dirac measure with mass at $x$ and  ${\cal M}(X)$
 stands for the set of all positive
 Radon  measures  on
${\cal B}(X)$. The space $\Gamma$ can be endowed with the
relative topology as a subset of the space ${\cal M}(X)$ with the
vague topology, i.e., the weakest topology on $\Gamma$ with
respect to which  all maps
$$\Gamma\ni\gamma\mapsto\la\varphi,\gamma\ra:=\int_X\varphi(x)\,\gamma(dx)
=\sum_{x\in\gamma}\varphi(x),\qquad\varphi\in C_0(X),$$ are continuous. Here, $C_0(X)$ denotes the set of all continuous functions
on $X$ with compact support.
We shall denote  the Borel $\sigma$-algebra on $\Gamma$ by
${\cal B}(\Gamma)$. If $\Xi$ is a subset of $\Gamma$, we shall denote by ${\cal B}(\Xi)$
the trace $\sigma$-algebra of ${\cal B}(\Gamma)$ on $\Xi$.

Let us fix any $x_0\in X$ and denote by $B(r):=B(x_0,r)$ the closed ball in $X$ of radius $r>0$, centered at $x_0$. For each $\alpha\ge1$, we define
$$\Theta_\alpha:= \big\{\,\gamma\in\Gamma\mid \exists K\in\N:\,
\forall r\in\N:\ |\gamma_{B(r)}|\le K\operatorname{vol}(B(r))^\alpha
 \,\big\},$$ where $\operatorname{vol}(B(r))$ denotes the volume of $B(r)$.
 We evidently have that $\Theta_{\alpha_1}\subset\Theta_{\alpha_2}$ if $\alpha_2\ge\alpha_1\ge1$. Denote also
 $$\Theta:=\bigcup_{\alpha\ge1}\Theta_\alpha.$$
  It easy to see that $\Theta_\alpha\in{\cal B}(\Gamma)$ for each $\alpha\ge1$, hence $\Theta\in{\cal B}(\Gamma)$.

For each $z>0$, let $\pi_z$ denote the Poisson measure on $(\Gamma,{\cal
B}(\Gamma))$ with intensity measure $z\,dx$. This measure can be
characterized by its Laplace transform
\begin{equation}\label{ewrwerewrwe5}\int_{\Gamma}
e^{\la\varphi,\gamma\ra}\,\pi_z(d\gamma)
=\exp\bigg(\int_X(e^{\varphi(x)}-1)\,z\,dx\bigg),\qquad
\varphi\in C_0(X).\end{equation} We refer to e.g.\ \cite{Kingman}
for
a detailed discussion of the construction of the Poisson measure
on the configuration space. By e.g.\ \cite{NZ}, we have, for each $z>0$, $$\pi_z(\Theta_1)=1$$ (in fact, by \cite{NZ},
 the Poisson measure $\pi_z$ is concentrated on those configurations $\gamma\in\Gamma$ for which  $\lim_{r\to\infty}|\gamma_{B(r)}|/\operatorname{vol}(B(r))=z$).

Using multiple stochastic integrals with respect to the Poisson random measure, one constructs a unitary isomorphism $$I_z:{\cal F}(L^2(X,z\,dx))
\to L^2(\Gamma,\pi_z),$$ see e.g.\ \cite{Surgailis1}.
Here, ${\cal F}(L^2(X,z\,dx))$ denotes the symmetric Fock space
over $L^2(X,z\,dx)$, i.e., $$ {\cal F}(L^2(X,z\,dx))=\bigoplus_{n=0}
^\infty {\cal F}_n(L^2(X,z\,dx)),$$ where $$ {\cal F}_n(L^2(X,x\,dx))
:= L^2(X,z\,dx)^{\odot n}n!,$$ $\odot$ standing for symmetric tensor product.

Let us recall the following  result of Surgailis \cite{Surgailis1}.
Let $A$ be a contraction in $L^2(X,z\,dx)$, and let $\operatorname{Exp}(A)$ denote the second quantization of $A$, i.e., $\operatorname{Exp}(A)$ is the contraction in ${\cal F}(L^2(X,z\,dx))$ given by
\begin{align*}
&\operatorname{Exp}(A)\restriction {\cal F}_0(L^2(X,z\,dx))=\pmb1,\\
&\operatorname{Exp}(A) f^{\otimes n}=(Af)^{\otimes n},\quad f\in L^2(X,z\,dx),\ n\in\N.
\end{align*} We shall keep the notation $\operatorname{Exp}(A)$ for the image of this operator under the isomorphism $I_z$.

In the following, we shall restrict out attention to the case of a self-adjoint $A$ (though the general case may be treated by an easy modification of the results below).

\begin{th}[Surgailis \cite{Surgailis1}] Let $A$ be a self-adjoint contraction in $L^2(X,z\,dx)$\rom. Then the operator $\operatorname{Exp}(A) $ is
positivity preserving if and only if $A$ is sub-Markov\rom. In the latter case\rom, $\operatorname{Exp}(A) $ is Markov\rom.

\label{jghyiu}
\end{th}

We recall that $A$ being sub-Markov means that $0\le Af \le1$ a.e.\ for each $0\le f\le1$ a.e., $f\in L^2(X,z\,dx)$. If, additionally, $Af_n \nearrow 1$ a.e.\ for some sequence $f_n\nearrow1$, $f_n\in L^2(X,z\,dx)$, then $A$ is called Markov.

Let $(T_t)_{t\ge0}$ be a Markov semigroup in $L^2(X,z\, dx)$, and let $p_t(x,\cdot)$, $t\ge0$, $x\in X$, be a corresponding Markov semigroup of kernels. Consider the semigroup $(\operatorname{Exp}(T_t))_{t\ge0}$ in $L^2(\Gamma,\pi_z)$,  which is Markov by Theorem~\ref{jghyiu}. Note that each operator $\operatorname{Exp}(T_t)$ is defined only $\pi_z$-{\it almost everywhere}. We are now interested in an explicit {\it point-wise} realization of a Markov semigroup of kernels which would correspond to the semigroup $(\operatorname{Exp}(T_t))_{t\ge0}$.

We consider the infinite product
$X^\N$
with the cylinder $\sigma$-algebra on it, denoted by ${\cal C}(X^\N)$. Let us recall the construction of a probability measure on $\Gamma$
through a product measure on $X^\N$, see \cite{VGG,KLR}.

We define ${\cal A}\in{\cal C}(X^\N)$ as the set of all elements $(x_n)_{n=1}^\infty\in X^\N$ such
that  $x_i\ne x_j$ when $i\ne j$, and  the sequence
$\{x_n\}_{n=1}^\infty$
 has no accumulation points in $X$.
 Let $$ D:=\{(x,y)\in X^2:\ x=y\}.$$
  Let $\nu_n$, $n\in\N$, be   probability measures on
$(X,{\cal B}(X))$ such that $$ \nu_n\otimes \nu_m(D)=0,\qquad n\ne m.$$
Consider the  product measure
$\nu:=\bigotimes_{n=1}^\infty\nu_n$  on $(X^\N,{\cal C}(X^\N))$.
Then, by the Borel--Cantelli lemma,  $\nu({\cal A})=1$
if and only if,  for each $r\in\N$, $$\sum_{n=1}^\infty \nu_n(B(r))<\infty.$$ In the latter case, we can consider $\nu$ as a probability measure on $({\cal A}, {\cal C}({\cal A}))$, where ${\cal C}({\cal A})$ denotes the trace $\sigma$-algebra of ${\cal C}(X^\N)$ on
${\cal A}$.  Define the mapping
\begin{equation}\label{ftyf}
{\cal A}\ni (x_n)_{n=1}^\infty\mapsto {\cal E}((x_n)_{n=1}^\infty):=
\sum_{n=1}^\infty \varepsilon_{x_n}\in\Gamma,
\end{equation}
which is measurable. Thus, we can define a probability measure $\mu$ on $(\Gamma,{\cal B}(\Gamma))$ as the image of $\nu$ under
the mapping \eqref{ftyf}. Evidently, the measure $\mu$ is independent
of the order of the $\nu_n$'s.

Assume now that \begin{equation}\label{hguyf} p_t(x,\cdot)\otimes p_t(y,\cdot)(D)=0,\qquad x\ne y, \ t>0,\end{equation}
and let $\gamma\in\Gamma$ be such that
\begin{equation}\label{guygiier} \sum_{x\in\gamma}p_t(x,B(r))<\infty,\qquad t>0,\ r\in\N.\end{equation} Then, for each $t\ge0$,
we define ${\bf P}_{t}(\gamma,\cdot)$ as the probability measure on $\Gamma$ given through the product measure $${\Bbb P}_{t}((x_n)_{n=1}^\infty,\cdot):=\bigotimes_{n=1}^\infty p_t(x_n,\cdot),$$ where $\{x_n\}_{n=1}^\infty$ is an arbitrary numeration of the elements of $\gamma$.

In what follows, we will always assume that the manifold $X$
satisfies the following condition: there exist $m\in\N$ and $C>0$ such that \begin{equation}\label{gyufr}
\vol(B(\beta r)) \le C \beta^m\vol(B(r)),\qquad r>0,\ \beta\ge1.
\end{equation}
By e.g.\ \cite[Proposition 5.5.1]{Davies}, if $X$
has non-negative Ricci curvature, then \eqref{gyufr} is satisfied with $C=1$ and $m$ being equal to the dimension of $X$.

\begin{th}\label{kjfvtft}  Let $(p_t)_{t\ge0}$ be a Markov semigroup of kernels on $X$ satisfying \eqref{hguyf} and let $\alpha\ge1$\rom.
 Assume that
\begin{equation}\label{gvjuvf}
\exists \epsilon>0:\ \forall t\in(0,\epsilon)\, \forall \delta>0:\quad \sum_{n=1}^\infty \sup_{x\in X}p_t(x,B(x,\delta n^{1/(\alpha m)})^c)<\infty.
\end{equation}
Then\rom, each $\gamma\in\Theta_\alpha $ satisfies \eqref{guygiier}\rom, so that ${\bf P}_{t}(\gamma,\cdot)$ is a probability measure on $\Gamma$ for each $t\ge0$\rom, and furthermore\rom,  $({\bf P}_t)_{t\ge0}$ is a Markov semigroup of kernels on $(\Theta_\alpha,{\cal B}(\Theta_\alpha))$\rom.

Additionally\rom, for each $z>0$ and $t>0$  and $F\in L^2(\Gamma,\pi_z)$\rom,
the function $$ \Theta_\alpha\ni\gamma\mapsto \int_{\Theta_\alpha}F(\xi){\bf P}_t(\gamma,d\xi)$$ is a $\pi_z$-version
of the function $\operatorname{Exp}(T_t)F\in L^2(\Gamma,\pi_z)$.

\end{th}

\begin{rem}\label{vfgtfuut}{\rm
Under the assumptions of Theorem \ref{kjfvtft}, if the condition \eqref{gvjuvf} is satisfied for all $\alpha\ge1$, then
$({\bf P}_t)_{t\ge0}$  becomes a Markov semigroup of kernels on $(\Theta,{\cal B}(\Theta))$.

}\end{rem}

\noindent{\it Proof of Theorem\/} \ref{kjfvtft}.
 For any $x\in X$, denote $|x|:=\operatorname{dist}(x_0,x)$, where $\operatorname{dist}$ denotes the Riemannian distance.

\begin{lem}\label{ysery} Let the conditions of Theorem \rom{\ref{kjfvtft}} be satisfied\rom. Then\rom, for any $\gamma\in\Theta_\alpha$ and $t\in(0,\epsilon)$, $$\sum_{x\in\gamma} p_t(x,B(x,|x|/2)^c)<\infty.$$
\end{lem}

\noindent{\it Proof}.   Fix any $\gamma\in\Theta_\alpha$ and choose any
numeration $\{x_n\}_{n=1}^\infty$ of points of $\gamma$  such that
$|x_{n+1}|\ge |x_n|$, $n\in\N$. Since $\gamma\in\Theta_\alpha$, there exists
$K\in\N$ for which \begin{equation}\label{hg}|\gamma_{B(r)}|\le
K\vol(B(r))^\alpha,\qquad r\in\N.\end{equation} Define
\begin{equation}\label{gfuf} r(n):=\max\{i\in\Z_+:\
i<(n/(KC^\alpha\vol(B(1))^\alpha))^{1/(\alpha m)}\},\qquad
n\in\N,\end{equation} where $C$ is the constant from \eqref{gyufr}. By
\eqref{gyufr}, \eqref{hg}, and \eqref{gfuf},  \begin{align*}
|\gamma_{B(r(n))}|&\le K\vol(B(r(n)))^\alpha\\
&\le KC^\alpha r(n)^{\alpha m} \vol(B(1))^\alpha\\ &<n,\quad n\in\N.\end{align*}
Therefore, $x_n\not\in B(r(n))$. Hence, by \eqref{gfuf},
$$|x_n|>r(n)>(n/(KC^\alpha\operatorname{vol}(B(1))^\alpha))^{1/(\alpha m)}-1. $$
Therefore, to prove the lemma, it suffices to show that
$$ \sum_{n=1}^\infty p_t(x_n,B(x_n, ((n/KC^\alpha\vol(B(1))^\alpha)^{1/(\alpha m)}-1)/2)^c)<\infty, $$
 which evidently  follows
 from \eqref{gvjuvf}.\quad $\square$

 Fix any $(x_n)_{n=1}^\infty\in{\cal E}^{-1}(\Theta_\alpha)$, where the mapping ${\cal E}$ is given by \eqref{ftyf}.
 Denote
 \begin{equation}\label{gyuuy} {\cal A}_n:=\big\{\,(y_k)_{k=1}^\infty \in X^\N:\ y_n\in B(x_n,|x_n|/2)\,\big\},\qquad n\in\N.\end{equation}  By Lemma \ref{ysery}, \eqref{gyuuy} and the Borel--Cantelli lemma,
 \begin{equation}\label{ghfd}
{\Bbb P}_t \left((x_n)_{n=1}^\infty,\liminf_{n\to\infty}{\cal A}_n\right)=1 ,\qquad \ t\in(0,\epsilon).
 \end{equation}
  Next, define $$ {\cal A}':= \liminf_{n\to\infty}{\cal A}_n\cap \{(y_n)_{n=1}^\infty\in X^\N:\ y_i\ne y_j\text{ if }i\ne j,\  i,j\in\N)\}.$$
 By \eqref{hguyf} and \eqref {ghfd},
  \begin{equation}\label{gfdyggy}
{\Bbb P}_t ((x_n)_{n=1}^\infty, {\cal A}')=1,\qquad \ t\in(0,\epsilon).
 \end{equation}
 We evidently have $$ |x_n|\to\infty\quad \text{as }n\to\infty,$$
   and therefore ${\cal A}'\subset{\cal A}$. Hence, condition \eqref{guygiier}  is satisfied for $\gamma={\cal E}((x_n)_{n=1}^\infty)$.

  Let us show that \begin{equation}\label{lhiouu}
  {\cal E}^{-1}(\Theta_\alpha)\subset{\cal A}'.\end{equation}
Indeed, fix any $(y_n)_{n=1}^\infty\in{\cal A}'$ and define $k\in{\Z}_+$ as the number of those $y_n$'s which do not belong to $B(x_n,|x_n|/2)$.   Then \begin{equation}\label{ufrtyu}
|{\cal E}((y_n)_{n=1}^\infty)_{B(r)}|\le |{\cal E}((x_n)_{n=1}^\infty)_{B(2r)}|+k.
\end{equation}
Then, by   \eqref{gyufr} and \eqref{ufrtyu}, we have, for each $r\in\N$,\begin{align*}
|{\cal E}((y_n)_{n=1}^\infty)_{B(r)}|&\le K\vol(B(2r))^\alpha+k\\
&\le K(C2^m\vol(B(r)))^\alpha+k\\ &\le K'\vol(B(r))^\alpha
\end{align*}
  for some $K'\in\N$ which is idependent of $r$
(note that $\vol(B(r))\to\infty $ as $r\to\infty$ since $X$ has infinite volume).   Hence ${\cal E}((y_n)_{n=1}^\infty)\in\Theta_\alpha$.

  Thus, by \eqref{gfdyggy} and \eqref{lhiouu},
  \begin{equation}\label{yrtsrea}
 {\Bbb P}_t((x_n)_{n=1}^\infty,{\cal E}^{-1}(\Theta_\alpha))=1,\qquad (x_n)_{n=1}^\infty\in{\cal E}^{-1}(\Theta_\alpha),\ t\in(0,\epsilon).
  \end{equation}
   Then, it easily follows from \eqref {yrtsrea}
and the construction of $ {\Bbb P}_t((x_n)_{n=0}^\infty,\cdot)$ that $({\Bbb P}_t)_{t\ge0}$ is a Markov semigroup of kernels on
${\cal E}^{-1}(\Theta_\alpha)$. Therefore, $({\bf P}_t)_{t\ge0}$ is a Markov semigroup of kernels on $\Theta_\alpha$.

 The proof of the last statement of the theorem is quite analogous to the proof of \cite[Theorem~5.1]{KLR}, so we only outline it. Let \begin{equation}\label{iguyfr} {\cal D}:=\{\varphi\in C_0(X): -1<\varphi\le0\}.\end{equation} Then, for any $\varphi\in {\cal D}$, $\gamma\in\Theta_\alpha$, and $t\ge0$, we easily get from the definition of ${\bf P}_t(\gamma,\cdot)$:
\begin{align}&\int_{\Theta_\alpha} \exp[\langle \log(1+\varphi),\xi\rangle]\,{\bf P}_t(\gamma,d\xi)\notag\\
&\qquad =\prod_{x\in\gamma}\int_X(1+\varphi(y))\,p_t(x,dy)\notag\\
&\qquad =\prod_{x\in\gamma} (1+(T_t\varphi)(x))\notag\\
&\qquad =\exp[\langle\log(1+T_t\varphi),\gamma\rangle].\label{drdrdr}
\end{align}

  Next, it is well known (see e.g.\ \cite[Corollary~2.1]{Surgailis1}) that, for any $\varphi\in{\cal D}$,
  \begin{equation}\label{iyufu}
  I_z^{-1}\exp[\langle\log(1+\varphi),\cdot\rangle]=\exp\bigg[\int_X \varphi(x)\,z\,dx\bigg]
  \big((1/n!)\varphi^{\otimes n}\big)_{n=0}^\infty.
  \end{equation}
  Since $ \int_X (T_t\varphi )(x)\,z\,dx=\int_X \varphi(x)\,z\,dx$, it follows from \eqref{drdrdr}, \eqref{iyufu} and the definition of $\operatorname{Exp}(T_t)$ that the statement   holds for $F=\exp[\langle\log(1+\varphi),\cdot\rangle]$, $\varphi\in {\cal D}$. From here, analogously to the proof of \cite[Theorem~5.1]{KLR}, we conclude the statement  in the general case.\quad $\square$
   \vspace{2mm}

   Let us now outline the case where the semigroup $(T_t)_{t\ge0}$ is sub-Markov, but not Markov. We shall assume for simplicity that \begin{equation}\label{igtr}
   \inf_{x\in X}p_t(x,X)>0,\qquad t>0.   \end{equation}
     Let $\hat X:= X\cup\{\Delta\}$ be a one-point extension of $X$, and,
     as usual, consider $(\hat p_t)_{t\ge0}$ as the extension of $(p_t)_{t\ge0}$ to a Markov semigroup of kernels on $\hat X$.
    Let conditions \eqref{hguyf} and  \eqref{gvjuvf} be satisfied. Consider the mapping
    $$\hat X^\N\ni(x_n)_{n=1}^\infty\mapsto \hat{\cal E}((x_n)_{n=1}^\infty):=\sum_{n=1}^\infty \pmb 1_X
      (x_n)\varepsilon_{x_n}.$$ Then  $\hat {\Bbb P}_t((x_n)_{n=1}^\infty,\cdot):=\bigotimes_{n=1}^\infty \hat p_t(x_n,\cdot)$,
      $t\ge0$, is a Markov semigroup of kernels on $\hat{\cal E}^{-1}(\Theta_\alpha)$ (notice that condition \eqref{igtr} guarantees that,
     for $\hat{\Bbb P}_t((x_n)_{n=1}^\infty,\cdot)$-a.e.\ $(y_n)_{n=1}^\infty\in\hat X^\N$, an infinite number
     of $y_n$'s belong to $X$).

     Set $\hat{\bf P}_t(\gamma,\cdot)$ to be the image of $\hat {\Bbb P}_t((x_n)_{n=1}^\infty,\cdot)$ under $\hat{\cal E}$, where $\gamma=\{x_n\}_{n=1}^\infty\in\Theta_\alpha$,
      $t\ge0$. Then, $(\hat{\bf P}_t)_{t\ge0}$ becomes a Markov semigroup of kernels on $\Theta_\alpha$.

  Next, we denote by $\pi_{z,t}$ the Poisson random measure over $X$ with intensity measure $(1-p_t(x,X))\,z\, dx$
 (notice that $\pi_{z,t}$ is concentrated on finite or infinite configurations in $X$, depending on whether the integral $\int_X(1-p_t(x,X))\,z\, dx$ is finite or infinite.) Define ${\bf P}_{z,t}(\gamma,\cdot)$ as the convolution of the measures $\hat {\bf P}_t(\gamma,\cdot)$ and $\pi_{z,t}$, i.e.,
 \begin{equation}\label{jhgftf} {\bf P}_{z,t}(\gamma,A)=\int \hat {\bf P}_t(\gamma,d\xi_1)\int \pi_{z,t}(d\xi_2)
 \pmb1_A(\xi_1+\xi_2).\end{equation}
 Note that ${\bf P}_{z,t}(\gamma,\cdot)$ is indeed concentrated on $\Gamma$, since for any fixed $\xi_1\in\Gamma$, the $\pi_{z,t}$
 probability of those $\xi_2$ which satisfy
$\xi_1\cap\xi_2\not=\varnothing$  is equal to zero. Furthermore, we have ${\bf P}_{z,t}(\gamma,\Theta_\alpha)=1$ for each $t\ge0$
and $\gamma\in\Theta_\alpha$. Indeed, for each $t>0$, $\pi_{z,t}$ is either concentrated on finite configurations or  $\pi_{z,t}(\Theta_\alpha)=1$, the latter being a consequence of the estimate $1-p_t(x,X)\le1$ and  the support property of a Poisson measure \cite{NZ}. Now, the equality  ${\bf P}_{z,t}(\gamma,\Theta_\alpha)=1$ follows from the definition \eqref{jhgftf}.

 Next, for any $\varphi\in {\cal D}$ (see \eqref{iguyfr}), $\gamma\in\Theta_\alpha$, and $t\ge0$, we have
  \begin{align}
  &\int_{\Theta_\alpha}\exp[\log(1+\varphi),\xi\rangle]
\, {\bf P}_{z,t}(\gamma,d\xi)\notag\\  &\quad =\int\exp[\langle \log(1+\varphi),\xi_1\rangle]\,\hat {\bf P}_t(\gamma,d\xi)\int\exp[\log(1+\varphi),\xi_2\rangle]\,\pi_{z,t}(d\xi_2)\notag\\
&\quad =\bigg(\prod_{x\in\gamma}\bigg(1-p_t(x,X)+\int_X(1+\varphi(y))p_t(x,dy)\bigg)\bigg)\exp\bigg[
\int_X\varphi(x)(1-p_t(x,X))\,z\,dx
\bigg]\notag\\ &\qquad= \exp[\langle\log(1+T_t\varphi),\gamma
\rangle] \exp\bigg[\int_X(\varphi(x)-(T_t\varphi)(x))\,z\,dx\bigg].\label{fytrde}
  \end{align}
    From here we conclude that \ $\Theta_\alpha\ni\gamma\mapsto\int_{\Theta_\alpha}F(\xi)\,{\bf P}_{z,t}(\gamma,d\xi)$
    is a $\pi_z$-version of $\operatorname{Exp}(T_t)F$.

Finally,  using \eqref{fytrde}, we have, for any $\varphi\in {\cal D}$, $t,s>0$, and $\gamma\in\Theta_\alpha$,
 $$\int_{\Theta_\alpha}\int_{\Theta_\alpha} \exp[\langle\log(1+\varphi),\xi_2\rangle]\,{\bf P}_{z,s}(\xi_1,d\xi_2)\,{\bf P}_{z,t}(\gamma,d\xi_1)=\int_\Theta\exp[\langle\log(1+\varphi),\xi\rangle]\,{\bf P}_{z,s+t}(\gamma,d\xi),$$
from where it easily follows that $({\bf P}_{z,t})_{t\ge0}$ is a Markov semigroup of kernels on $\Theta_\alpha$.

   \begin{rem}\label{ftuyfrut}\rom{ Note that, in the case where the semigroup $(T_t)_{t\ge0}$ is Markov, the construction
   of the Markov semigroup of kernels $({\bf P}_t)_{t\ge0}$ is independent of $z$, whereas in the case where $(T_t)_{t\ge0}$
   is sub-Markov, the $({\bf P}_{z,t})_{t\ge0}$ does depend on $z$. }

   \end{rem}

   \section{Non-equilibrium independent infinite particle process}\label{hgig}

    Our next aim is to study the Markov process corresponding to the Markov semigroup of kernels $({\bf P}_t)_{t\ge0}$,
    respectively $({\bf P}_{z,t})_{t\ge0}$.

    For a metric space $E$, we denote by $D([0,\infty),E)$ the space of all cadlag functions from $[0,\infty)$ to $E$, i.e., right continuous functions on $[0,\infty)$ having left limits on $(0,\infty)$. We  equip $D([0,\infty),E)$  with the cylinder $\sigma$-algebra ${\cal C}(D([0,\infty),E))$ constructed through the Borel $\sigma$-algebra ${\cal B}(E)$.

    In what follows, we will  assume that a Markov process on $X$ corressponding to the semigroup $(T_t)_{\ge0}$ has cadlag paths.  The latter, in particular, holds if the kernels
    $(p_t)_{t\ge0}$ determine a Feller semigroup (see e.g.\ \cite[Chapter~2, Theorem~2.7]{EK}).

   We first consider the case where the semigroup $(T_t)_{t\ge0}$ is Markov. For each $x\in X$, let $P^x$ denote the distribution of the Markov process $(X_t)_{t\ge0}$ corresponding to $(T_t)_{t\ge0}$ which starts at $x$. By our assumption, each $P^x$ is a probability measure on  $D([0,\infty),X)$.    We will also assume that
   \begin{equation}\label{uytd}
   P^x\otimes P^y((X_t^{(1)},X_t^{(2)})_{t\ge0}\mid \exists t>0:\ X_t^{(1)}=X^{(2)}_t)=0,\qquad x\ne y,
   \end{equation}
  i.e., two independent Markov process starting at $x$ and $y$, $x\ne y$, will a.s.\ never meet. (If this condition is not satisfied, all the results below remain true, but for a corresponding space of multiple configurations.) Notice that
  condition \eqref{uytd} is stronger than \eqref{hguyf}.

    For each $x\in X$ and $r>0$, denote by $\tau_{B(x,r)^c}$ the hitting time of $B(x,r)^c$:
    $$ D([0,\infty),X)\ni\omega\mapsto \tau_{B(x,r)^c}(\omega):=
    \inf\{t>0:\, \omega(t)\in B(x,r)^c\}.$$ Let $\alpha\ge1$ and assume:
    \begin{equation}\label{ufu}\exists \epsilon>0:\ \forall\delta>0:\
    \sum_{n=1}^\infty \sup_{x\in X} P^x(
    \tau_{B(x,\delta n^{1/(\alpha m)})^c}>\epsilon
    )<\infty.
        \end{equation} Condition \eqref{ufu} is evidently stronger than \eqref{gvjuvf}.

     Consider the space $D([0,\infty),X)^\N$ equipped with the cylinder $\sigma$-algebra  \linebreak ${\cal C}(D([0,\infty),X)^\N)$.
       Denote by $\Omega_{\alpha,1}$ the set of those $(\omega_n)_{n=1}^\infty\in D([0,\infty),X)^\N$ which satisfy
    the following conditions:
    \begin{itemize}
    \item[(i)]
   for all $t\ge0$, $\omega_i(t)\ne\omega_j(t)$ if $i\ne j$;
   \item[(ii)] $\{\omega_n(0)\}_{n=1}^\infty\in\Theta_\alpha$;
  \item[(iii)] there are only a finite number of $\omega_n$'s for which $\tau_{B(x_n,|x_n|/2)^c}(\omega_n)\le \epsilon$, where $\epsilon$ is as in \eqref{ufu}
  \end{itemize}
    It is easy to see that $\Omega_{\alpha,1}\in{\cal C}(D([0,\infty),X)^\N)$.

        Fix any $(x_n)_{n=1}^\infty\in{\cal E}^{-1}(\Theta_\alpha)$ and consider the product measure
    $${\Bbb P}^{(x_n)_{n=1}^\infty}:=\bigotimes_{n=1}^\infty P^{x_n}$$ on  $D([0,\infty),X)^\N$.
        Absolutely analogously to the proof of Theorem \ref{kjfvtft} we conclude from \eqref{uytd} and \eqref{ufu} that
    \begin{equation}\label{ufyfyufy}
    {\Bbb P}^{(x_n)_{n=1}^\infty}(\Omega_{\alpha,1})=1.\end{equation}
    Furthermore, it follows from  the proof of Theorem \ref{kjfvtft} that, for each $(\omega_n)_{n=1}^\infty\in\Omega_{\alpha,1}$, we have $\{\omega_n(t)\}_{n=1}^\infty\in\Theta_\alpha$ for all $t\in(0,\epsilon]$.

    For each $k\in\N$, we now recurrently define  $\Omega_{\alpha,k+1}$ as the subset of $\Omega_{\alpha,k}$ consisting of those $(\omega_n)_{n=1}^\infty$ which satisfy the following condition: there are only a finite number of $\omega_n$'s for which
    $$ \tau_{B(\omega_n(\epsilon k),|\omega_n(\epsilon k)|/2)^c}(\omega_n(\epsilon k+\cdot))\le\epsilon.$$
 Since under $P^{x}$,  $X_t(\omega)=\omega(t)$, $t\ge0$,  is a time homogeneous Markov process starting at $x$,
 $x\in X$, we conclude, analogously to the above, that $$ {\Bbb P}^{(x_n)_{n=1}^\infty}(\Omega_{\alpha,k})=1,\qquad k\in\N.$$
 Therefore, $$ {\Bbb P}^{(x_n)_{n=1}^\infty}(\Omega_\alpha)=1, \quad\Omega_\alpha:=\bigcap_{k\in\N}\Omega_{\alpha,k}.$$
Hence, we can consider ${\Bbb P}^{(x_n)_{n=1}^\infty}$ as a probability measure on $\Omega_\alpha$
equipped with the trace $\sigma$-algebra of ${\cal C}(D([0,\infty),X^\N)$ on $\Omega_\alpha$.

Fix any  $(\omega_n)_{n=1}^\infty\in\Omega_\alpha$. Then, we evidently have $\{\omega_n(t)\}_{n=1}^\infty\in\Theta_\alpha$ for all $t\ge0$.
 Furthermore, for any compact $\Lambda\subset X$ and for any $T>0$, there are only a finite number of $\omega_n$'s which meet $\Lambda$ during the time interval $[0,T]$. Therefore, $\{\omega_n(\cdot)\}_{n=1}^\infty\in D([0,\infty),\Theta_\alpha)$, where $\Theta_\alpha$ is equipped with the relative topology as a subset of $\Gamma$. Thus, the following mapping is well-defined:
\begin{equation}\label{kjbgi} \Omega_\alpha\ni (\omega_n(\cdot))_{n=1}^\infty \mapsto \sum_{n=1}^\infty \varepsilon_{\omega_n(\cdot)} \in D([0,\infty),\Theta_\alpha).\end{equation}
Furthermore, it easily follows from the definition of a cylinder $\sigma$-algebra that the mapping \eqref{kjbgi} is measurable. Thus, we can consider the image of the measure  ${\Bbb P}^{(x_n)_{n=1}^\infty}$ under \eqref {kjbgi}. This probability measure on $D([0,\infty),\Theta_\alpha)$  will be denoted by ${\bf P}^{\{x_n\}_{n=1}^\infty}$ (note that this measure is, indeed, independent of the numeration of points of $\{x_n\}_{n=1}^\infty\in\Theta_\alpha$).

 Next, it is easy to see that the finite-dimensional distributions of
       ${\bf P}^\gamma$, $\gamma\in\Theta_\alpha$, are given through the Markov semigroup of kernels ${\bf P}_t(\gamma,\cdot)$ on $\Theta_\alpha$ (see Theorem~\ref{kjfvtft}). From here, analogously to \cite[Theorem~8.1]{KLR}, we get the following

\begin{th}\label{oihiuig} Let $(p_t)_{\ge0}$ be a Markov semigroup of kernels on $X$ and let  \eqref{uytd} and \eqref{ufu} hold.
Then there exists  a  time homogeneous Markov process
$${\bf M}=(\pmb\Omega,{\bf F},({\bf F}_t)_{t\ge0},(\pmb\theta_t)_{t\ge0},({\bf P}^\gamma)_{\gamma\in\Theta}, ({\bf X}_t)_{t\ge0}) $$
       on the state space $(\Theta_\alpha,{\cal B}(\Theta_\alpha))$ with cadlag paths and with  transition probability function $({\bf P}_t)_{t\ge0}$.

\end{th}

\begin{rem}
\rom{
In Theorem \ref {oihiuig}, ${\bf M}$ can be taken canonical, i.e., $\pmb\Omega =D([0,\infty),\Theta_\alpha)$, ${\bf X}(t)(\omega)=\omega(t)$ for  $t\ge0$ and  $\omega\in\pmb\Omega$, ${\bf F}_t=\sigma\{{\bf X}_s, 0\le s\le t\}$ for $t\ge0$, ${\bf F}=\sigma\{{\bf X}_t,t\ge0\}$, $(\pmb\theta_t\omega)(s)=\omega(s+t)$ for  $t,s\ge0$.
}\end{rem}

 \begin{rem}\rom{
 From the proof of Theorem \ref{oihiuig} we see that the Markov process $\bf M$ is a realization of the independent infinite particle process.  }\end{rem}

\begin{rem}\rom{
If the underlying Markov process on $X$ has continuous sample paths, then the Markov process ${\bf M}$ in Theorem \ref{oihiuig}
has continuos sample paths in $\Theta_\alpha$. }\end{rem}

     \begin{rem}\rom{
   Analogously to Remark  \ref{vfgtfuut}, we note that if, under the assumptions of Theorem~\ref{oihiuig}, condition \eqref{ufu}
   is satisfied for each $\alpha\ge1$, then  $\Theta_\alpha$ can be replace by $\Theta$
in the statement  of this theorem.

     }\end{rem}

 \begin{cor}\label{uyfrturf}
 The statement of Theorem \rom{\ref{oihiuig} } remains true if, instead of \eqref{ufu}, one demands that the following  stronger condition be satisfied:
 \begin{equation}\label{ufutrswer}\exists \epsilon>0:\ \forall\delta>0:\
    \sum_{n=1}^\infty
  \sup_{t\in(0,\epsilon]}\sup_{x\in X} p_t(x,
    B(x,\delta n^{1/(\alpha m)})^c  )<\infty.
        \end{equation}
 \end{cor}

\noindent {\it Proof}.  For any $x\in X$, $r>0$, and $\epsilon>0$, we have:
\begin{equation}\label{lif}
P^x(\tau_{B(x,r)^c}>\epsilon)\le 2\sup_{t\in(0,\epsilon]}\sup_{x\in X}
p_t(x,B(x,r/2)^c).
\end{equation}
This estimate  follows by a straightforward generalization of (the proof of) \cite[Appendix A, Lemma 4]{Nelson}
(see also \cite[Lemma 8.1]{KLR}) to the case of an arbitrary Markov process on $X$ with cadlag paths.
Now, by \eqref{lif}, condition \eqref{ufutrswer}
implies \eqref{ufu}. \quad $\square$

  Let us consider the case where $(p_t)_{t\ge 0}$ is sub-Markov.
  We will assume that, for each $x\in X$,
  $p_{t}(x,X)$ is continuously differentiable in $t\in[0,\infty)$ and
   there exists $\delta>0$ such that
   \begin{equation}\label{hjkg}
   \left| \frac{\partial}{\partial t}\, p_{t,x}(X)\right|\le C,\qquad t\in[0,\delta],\ x\in X,\end{equation} for some $C>0$.

  We set \begin{equation}\label{jkguyg}g(x):=-\frac{\partial}{\partial t}\,p_{t,x}(X)\Big|_{t=0},\qquad x\in X,\end{equation}
  which is a bounded non-negative function.
    Using the semigroup property of $(p_t)_{t\ge0}$, we then conclude from \eqref{hjkg} and \eqref{jkguyg} that \begin{equation}\label{ufytd}
  \frac{\partial}{\partial t}\,p_{t,x}(X)=-\int_X g(y) p_{t,x}(dy),\qquad t>0,\ x\in X.
  \end{equation}

  We will now assume that \eqref{igtr},  \eqref{uytd}, and \eqref{ufu}
  hold. Analogously to the Markovian case, for each  $\gamma\in\Theta_\alpha$, we define a probability measure $\hat{\bf P}^\gamma$ on $D([0,\infty),\Theta_\alpha)$ through the mapping \begin{equation}\label{ufyt}\hat\Omega_\alpha\ni(\omega_n(\cdot))_{n=1}^\infty \mapsto\sum_{n=1}^\infty
  \pmb 1_X(\omega_n(\cdot))\varepsilon_{\omega_n(\cdot)}\in D([0,\infty),\Theta_\alpha),\end{equation}
  where $\hat\Omega_\alpha$ is a corresponding subset of $D([0,\infty),\hat X)^{\N}$ (compare with \eqref{kjbgi}).

  Let $\Pi_z$ denote the Poisson random measure over $X\times[0,\infty)$ with intensity measure $g(x)\,z\, dx\,dt$. Since the measure $dx$ has no atoms in $X$, $\Pi_z$ is concentrated on the set of those configurations $\xi$ which satisfy the following condition: for any different $(x_1,t_1),(x_2,t_2)\in\xi$, we have $x_1\ne x_2$. Any such configuration can be represented as the disjoint union $$\xi=\bigcup_{k=1}^\infty\xi^{(k)},$$ where each $\xi^{(k)}$ is a configuration in $X\times [k-1,k)$, and to each $\xi^{(k)}$ there corresponds a configuration $\gamma^{(k)}$ in $X$ that is obtained by taking the $X$-components of the points from $\xi^{(k)}$.

Denote  $\overline\Theta_\alpha:=\Theta_\alpha\cup\Gamma_{\mathrm{fin}}$, where
$\Gamma_{\mathrm{fin}}$ denotes the set of all finite configurations in $X$. We endow $\overline\Theta_\alpha$ with the vague topology.  We note that, since the function $g(x)$ is bounded, we have $\Pi_z$-a.s.\ that $\gamma^{(k)}\in\overline\Theta_\alpha$.

For each $\xi^{(k)}$, we now construct a probability measure ${\bf M}^{\xi^{(k)}}$ on $D([0,\infty_,\overline\Theta_\alpha)$. This measure is defined in the same way as the measure $\hat{\bf P}^{\gamma^{(k)}}$, but only, instead of \eqref{ufyt}, one uses the mapping
$$ (\omega_n(\cdot))_{n\ge1}\mapsto\sum_{n\ge1}\pmb 1_{[t_n,\infty)}
(\cdot)\pmb1_X(\omega_n(\cdot-t_n))\varepsilon_{\omega_n(\cdot-t_n)},$$ where
$\xi^{(k)}=\{(x_n,t_n)\}_{n\ge1}$, so that $\gamma^{(k)}=\{x_n\}_{n\ge1}$.

For each $\gamma\in\Theta_\alpha$, we now define a probability measure ${\bf P}_z^\gamma$ on $D([0,\infty),\Theta_\alpha)$
by setting, for each $C\in{\cal C}(D([0,\infty),\Theta_\alpha)$,
$$ {\bf P}_z^\gamma(C):=\int\hat{\bf P}^\gamma(d\pmb\omega^{(0)})\int\Pi_z(d\xi)\int\bigg(\bigotimes_{k=1}^\infty
{\bf M}^{\xi^{(k)}}\bigg)(d\pmb\omega^{(1)}(\cdot),d\pmb\omega^{(2)}(\cdot),\dots)\pmb1_C\bigg(
\sum_{k=0}^\infty\pmb\omega^{(k)}(\cdot)
\bigg).$$
We  note that $\sum_{k=0}^\infty\pmb\omega^{(k)}(\cdot)$ indeed a.s.\ belongs to $D([0,\infty),\Theta_\alpha)$, since
for any $K\in\N$ and $t<K$, we have $\sum_{k=K}^\infty\pmb\omega^{(k)}(t)=0$.

Let us show that the finite-dimensional distributions of ${\bf P}_z^\gamma$ are given through the Markov semigroup of kernels
${\bf P}_{z,t}(\gamma,\cdot)$ on $\Theta_\alpha$.  For  $0=t_0<t_1<t_2<\dots< t_n$, $n\in\N$, and any $\varphi_1,\dots,\varphi_n\in {\cal D}$
(see \eqref{iguyfr}), we have
\begin{align}
&\int \Pi_z(d\xi) \int \bigg(\bigotimes_{k=1}^\infty {\bf M}^{\xi^{(k)}}\bigg)(d\pmb\omega^{(1)}(\cdot),d\pmb\omega^{(2)}(\cdot),
\dots) \prod_{i=1}^n\exp\bigg[\bigg\langle\log(1+\varphi_i),\sum_{k=1}^\infty \pmb\omega^{(k)}(t_i)\bigg\rangle\bigg]\notag\\
&\quad=\int\Pi_z(d\xi)\prod_{(x,\tau)\in\xi}\int\hat P^x(d\omega(\cdot))\prod_{i=1}^n\exp[\langle\log(1+\varphi_i),\pmb1_{[t,\infty)}
(t_i)\pmb1_X(\omega(t_i-\tau))\varepsilon_{\omega(t_i-\tau)}\rangle]\notag\\
&\quad=\exp\bigg[\int_X z\, dx\,g(x)\int_0^\infty d\tau\bigg(-1+\int\hat P^x(d\omega(\cdot))\notag\\ &\qquad\times
\prod_{i=1}^n \exp[\langle\log(1+\varphi_i),\pmb1_{[t,\infty)}
(t_i)\pmb1_X(\omega(t_i-\tau))\varepsilon_{\omega(t_i-\tau)}\rangle]\bigg)\bigg]\notag\\
&\quad=\exp\bigg[\sum_{j=1}^n \int_X z\,dx\,g(x)\int_{t_{j-1}}^{t_j} d\tau\bigg(-1+\int\hat P^x(d\omega(\cdot))\notag\\ &\qquad\times
\prod_{i=1}^n \exp[\langle\log(1+\varphi_i),\pmb1_X(\omega(t_i-\tau))\varepsilon_{\omega(t_i-\tau)}\rangle]\bigg)\bigg]\notag\\
&\quad=\exp\bigg[\sum_{j=1}^n \int_X z\,dx\,g(x)\int_{t_{j-1}}^{t_j} d\tau
\bigg(-1+(1-p_{t_j-\tau}(x,X))+\int_X p_{t_j-\tau}(x,dy)(1+\varphi_j(y))\notag\\
&\qquad\times\int_X\hat P^y(d\omega(\cdot))\prod_{i=j+1}^n
\exp[\langle\log(1+\varphi_i),\pmb1_X(\omega(t_i-t_j))\varepsilon_{\omega(t_i-t_j)}\rangle]\bigg)\bigg]\label{yiy}.\end{align}
Denote $$ F_j(x):=-1+(1+\varphi_j(x))\int\hat P^x(d\omega(\cdot))\prod_{i=j+1}^n \exp[\langle\log(1+\varphi_i),\pmb1_X(\omega(t_i-t_j))\varepsilon_{\omega(t_i-t_j)}\rangle],$$
where $x\in X$. Then we can proceed in \eqref{yiy} as follows:
\begin{align}
&=\exp\bigg[\sum_{j=1}^n \int_X z\,dx \, g(x) \int_{0}^{t_j-t_{j-1}} d\tau
\int_X p_{t_j-t_{j-1}-\tau}(x,dy)F_j(y)\notag\\ &=\exp\bigg[ \sum_{j=1}^n \int_0^{t_j-t_{j-1}}d\tau\int_X z\,dx\,F_j(x)\int p_{t_j-t_{j-1}-\tau}(x,dy)g(y)\bigg]\notag\\
&=\exp\bigg[\sum_{j=1}^n\int_X  F_j(x)(1-p_{t_j-t_{j-1}}(x,X))\,z\,dx\bigg],
\label{futyd}
\end{align} where we used \eqref{ufytd}.

  On the other hand, by \eqref{jhgftf},
  \begin{align}
&\int{\bf P}_{z,t_1}(\gamma,d\gamma_1)\int{\bf P}_{z,t_2-t_1}(\gamma_1,d\gamma_2)\dotsm\int{\bf P}_{z,t_n-t_{n-1}}(\gamma_{n-1},d\gamma_n)
\prod_{i=1}^n\exp[\langle\log(1+\varphi_i),\gamma_i\rangle]\notag\\
&\quad=\int\hat{\bf P}_{t_1}(\gamma,d\gamma_1)\int\hat{\bf P}_{t_2-t_1}(\gamma_1,d\gamma_2)\dotsm\int\hat{\bf P}_{t_n-t_{n-1}}(\gamma_{n-1},d\gamma_n)
\prod_{i=1}^n\exp[\langle\log(1+\varphi_i),\gamma_i\rangle]\notag\\
&\qquad\times \prod_{j=1}^n A_j,\label{ifgiufyg}
\end{align}
  where
  \begin{align}
  A_j&=\int\pi_{z,t_j-t_{j-1}}(d\theta_j)\exp[\langle\log(1+\varphi_j),\theta_j\rangle]\notag\\
&\quad\times
 \int\hat{\bf P}_{t_{j+1}-t_j}(\theta_j,d\eta_{j+1})\dotsm\int\hat{\bf P}_{t_n-t_{n-1}}(\eta_{n-1},d\eta_n)\prod_{k=j+1}^n \exp[\langle \log(1+\varphi_k),\eta_k\rangle]\notag\\
 &= \int\pi_{z,t_j-t_{j-1}}(d\theta_j) \prod_{x\in\theta_j}\bigg( (1+\varphi_j(x))\int\hat p_{t_{j+1}-t_j}(x,dy_1)\dotsm \int \hat p_{t_n-t_{n-1}}(y_{n-1},dy_n)\notag\\&\qquad\times \prod_{k=j+1}^n\exp[\langle \log(1+\varphi_k),
 \pmb1_X(y_k)\varepsilon_{y_k}\rangle]\bigg)\notag\\&=\exp\bigg[
 \int_X F_j(x)(1-p_{t_j-t_{j-1}}(x,X))\,z\,dx\bigg].  \label{rewawa}\end{align}
    By \eqref{yiy}--\eqref{rewawa}, we conclude that the finite-dimensional distributions of ${\bf P}_z^\gamma$ are indeed  given through the Markov semigroup of kernels
${\bf P}_{z,t}(\gamma,\cdot)$

  Thus, analogously to Theorem \ref{oihiuig}, we obtain a  time homogeneous Markov process ${\bf M}_z$ on the state space $\Theta_\alpha$ (or $\Theta$, provided \eqref{ufu} holds for all $\alpha\ge1$)  with cadlag paths and with  transition probability function $({\bf P}_{z,t})_{t\ge0}$.

\begin{rem}\rom{
Analogously to Remark \ref{ftuyfrut}, we see that, in the case where $(p_t)_{t\ge0}$ is Markov, the process ${\bf M}$
from Theorem \ref{oihiuig} has {\it any\/} Poisson measure $\pi_z$, $z>0$ as invariant measure, whereas, in the case of a sub-Markov $(p_t)_{t\ge0}$, only the measure $\pi_z$ is invariant for the process ${\bf M}_z$.
}\end{rem}

\section{Examples}\label{kugyug}

\subsection{Brownian motion on the configuration space}

Assume that $(T_t)_{t\ge0}$ is the heat semigroup  on $X$ with generator $\frac12\Delta^X$, the Laplace--Beltrami operator on $X$. We will denote by $p(t,x,y)$ the corresponding heat kernel on $X$ (see e.g.\ \cite{Davies}). We recall that the corresponding Markov process on $X$ is called Brownian motion on $X$.
In the case where the manifold $X$ is not stochastically complete, we will assume that condition \eqref{hjkg} is satisfied.

\begin{th}\label{ftydy}
Assume that the dimension of the manifold $X$ is $\ge2$.
Assume that $(p_t)_{t\ge0}$ is either Markov, or \eqref{igtr} and \eqref{hjkg} hold. Further assume that the heat kernel $p(t,x,y)$ satisfies the Gaussian upper bound for small values of $t$:
\begin{equation}\label{gfydytr} p(t,x,y)\le C t^{-n/2}\exp\bigg[-\frac{\operatorname{dist}(x,y)^2}{Dt}\bigg],\qquad t\in(0,\epsilon],\ x,y\in X,\end{equation} where $n\in\N$, $\epsilon>0$ and $C$ and $D$ are positive constants. Then the corresponding independent infinite particle process exists as a Markov process ${\bf M}^{\mathrm B}$ on $(\Theta,{\cal C}(\Theta))$ with either continuous paths if
$(p_t)_{t\ge0}$ is Markov, or cadlag paths if $(p_t)_{t\ge0}$ is
sub-Markov.
\end{th}

\begin{rem}\rom{
According to \cite{AKR3} and \cite{KLR}, the Markov process ${\bf M}^{\mathrm B}$
in Theorem \ref{ftydy} may be interpreted as a Brownian motion in  the configuration space over $X$.

Denote by ${\cal F}C_{\mathrm b}^\infty(C_0^\infty(X),\Theta)$ the set of all real-valued functions on $\Theta$ of the form
$F(\gamma)=g_F(\langle \varphi_1,\gamma\rangle,\dots, \langle \varphi_N,\gamma\rangle)$, where $g_F\in C_{\mathrm b}
^\infty(\R)$, $\varphi_1,\dots,\varphi_N\in C_0^\infty(X)$, $N\in\N$.  Assume first that $(p_t)_{t\ge0}$ is Markov. Then, the $L^2$-generator of the process ${\bf M}^{\mathrm B}$ has the following representation on the set ${\cal F}C_{\mathrm b}^\infty(C_0^\infty(X),\Theta)$
(which is a core for this operator):
$$ ({\bf L}^{\mathrm B}F)(\gamma)=\frac12\sum_{x\in\gamma} \Delta_x^XF(\gamma),$$ where $$\Delta^X_xF(\gamma):=
\Delta_y^X F(\gamma\setminus\{y\}\cup\{x\})\big|_{y=x},$$ see \cite{AKR3,KLR} for details. In the case where $(p_t)_{t\ge0}$ is sub-Markov, one can analogously show that the  $L^2$-generator of  ${\bf M}^{\mathrm B}$ is given on the set ${\cal F}C_{\mathrm b}^\infty(C_0^\infty(X),\Theta)$ by
$$ ({\bf L}^{\mathrm B}F)(\gamma)=\frac12\sum_{x\in\gamma} \Delta_y^X(F(\gamma\setminus \{ x \} \cup \{ y \} )-F(\gamma\setminus \{ x \} )\big|_{y=x}+\int_X (F(\gamma\cup\{x\})-F(\gamma))g(x)\,dx.$$
}\end{rem}

\noindent{\it Proof of Theorem\/} \ref{ftydy}. Since the dimension of $X$ is
$\ge2$, \eqref{uytd} is now satisfied, see e.g.\ (8.29) in  \cite{KLR}. Next,
by \cite[Lemma 8.2]{KLR},  \eqref{gfydytr}  implies that there exists $C>0$
such that $$ \sup_{t\in(0,\epsilon]}\sup_{x\in X} p_t(x,B(x,r)^c)\le C
e^{-r},\qquad r>0.$$ It follows from here that \eqref{ufutrswer} is satisfied.
Now, the theorem follows from  Corollary~\ref{uyfrturf}. \quad $\square$

 \subsection{Free Glauber dynamics on the configuration space}\label{gyfhh}

 Let $a:X\to[0,\infty)$ be a  bounded measurable function.
 Consider a Markov process on $X$ with a finite life time that corresponds to the semigroup $(T_tf)(x)=e^{-a(x)t}f(x)$. Thus, the process stays at a starting point for some random time and then dies, with $$p_t(x,\{x\})=p_t(x,X)=e^{-a(x)t},\qquad x\in X,\, t>0.$$ The function $g$ is now equal to $a$. Note also that $p_t(x,B(x,r)^c)=0$ for each $r>0$. Thus, all the above conditions are evidently satisfied and the corresponding infinite particle process exists as a Markov process ${\bf M}^{\mathrm G}$ on $(\Theta,{\cal C}(\Theta))$ with cadlag
 paths. This process can be interpreted as a free Glauber dynamics on the configuration space, or a birth-and-death process in $X$, see \cite{BCC,KL,KLR1}
 for details.
The $L^2$-generator of the process ${\bf M}^{\mathrm G}$ is given on its core
${\cal F}C_{\mathrm b}^\infty(C_0^\infty(X),\Theta)$ by $$ ({\bf L}^{\mathrm
G}F)(\gamma)=\sum_{x\in\gamma}a(x)(F(\gamma\setminus\{x\})-F(\gamma)) +\int_X
a(x)(F(\gamma\cup\{x\}-F(\gamma))\,dx,$$ see \cite{KL,KLR1}.

  \subsection{Free Kawasaki dynamics on the configuration space}

  We now consider a Markov jump process on $X$. The generator of this process  has the following
  representation on the set of all bounded functions on $X$:
  \begin{equation}\label{jufty} (Lf)(x) =\int_X (f(y)-f(x))\varkappa(x,y)\,dy,\end{equation} and we assume that $$ \varkappa:X^2\to[0,\infty)$$
  is a measurable function satisfying \begin{equation}
\label{oihiu}
\varkappa(x,y)=\varkappa(y,x),\qquad x,y \in X,
\end{equation}
  and \begin{equation}\label{jff}
  \lambda:= \sup_{x\in X}\int_X\varkappa(x,y)\,dy\in(0,\infty).
  \end{equation}

  Following \cite{EK}, we can explicitly construct this Markov process as follows. For each $x_0\in X$, let $\{Y(k),\ k=0,1,2,\dots\}$  be a Markov chain in $X$ starting at $x_0$,  with transition function
  $$\mu(x,dy)=\bigg(1-\frac1\lambda\,\int_X\varkappa(x,y)\,dy\bigg)\varepsilon_x(dy)+\frac1\lambda\,\varkappa(x,y)\,dy.$$
  Let $(Z_t)_{t\ge0}$ be an independent Poisson process with parameter $\lambda$. We now define  the Markov process
  $(X_t)_{t\ge0}$ starting at $x_0$ by $$ X_t:= Y(Z_t),\qquad t\ge0.$$
  By \cite{EK}, this process has generator \eqref{jufty}.

   \begin{th}\label{vgfuft}
  Assume that \eqref{oihiu} and \eqref{jff} hold and assume that there exist $C>0$ and $\alpha>m$ \rom($m$ being as in \eqref{gyufr}\rom) such that
  \begin{equation}\label{jufuyt}
  \sup_{x\in X}\int_{B(x,r)^c}\varkappa(x,y)\,dy\le \frac{C}{r^\alpha},\qquad r>0.
  \end{equation} Then the corresponding independent
   infinite particle process exists as a Markov process
   ${\bf M}^{\mathrm K}$ on $(\Theta_\alpha,{\cal B}(\Theta_\alpha))$ with cadlag paths.

\end{th}

  \begin{rem}
\rom{  According to \cite{KLR1}, the Markov process ${\bf M}^{\mathrm K}$ can be interpreted as a free Kawasaki dynamics on the configuration space. The $L^2$-generator of the process ${\bf M}^K$  has the following representation
on ${\cal F}C_{\mathrm b}^\infty(C_0^\infty(X),\Theta)$:
$$ ({\bf L}^{\mathrm K}F)(\gamma)=\int_X\gamma(dx)\int_X dy\, \varkappa(x,y)(F(\gamma\setminus\{x\}\cup\{y\})-F(\gamma)),$$
  see \cite{KLR1}.
  }\end{rem}

  \noindent {\it Proof}.
   By construction,  the process $(X_t)_{t\ge0}$ has cadlag paths in $X$. Furthermore, as easily seen, condition
  \eqref{uytd} is now satisfied.
  In a standard way, the process $(X_t)_{t\ge}$  leads to a Markov semigroup $(T_t)_{t\ge0}$ in $L^2(X,dx)$. Using \eqref{oihiu}
  and the construction of the process, we see that each $T_t$ is self-adjoint.

 By Theorem \ref{oihiuig}, it suffices to prove that condition \eqref{ufu} is satisfied. By using \eqref{jufuyt} and the construction of the process
 $(X_t)_{t\ge0}$,  we get, for any $\delta>0$.
  \begin{align*}
 &\sum_{n=1}^\infty \sup_{x\in X}P^x(\tau_{B(x,\delta n^{1/m})^c}>1)\\
 &\qquad =
 \sum_{n=1}^\infty\sup_{x\in X}\sum_{k=1}^\infty P^x(\tau_{B(x,\delta n^{1/m})^c}>1,\, Z_1=k)\\
 &\qquad \le \sum_{n=1}^\infty \sup_{x\in X}\sum_{k=1}^\infty P^x
 (Z_1=k,\, \exists i\in\{1,\dots,k\}:\, \operatorname{dist}(Y(i-1),Y(i))\ge \delta n^{1/m}/k)\\ &\qquad= \sum_{n=1}^\infty
 \sup_{x\in X}\sum_{k=1}^\infty
 e^{-\lambda}\frac{\lambda^k}{k!}P^x( \exists i\in\{1,\dots,k\}:\, \operatorname{dist}(Y(i-1),Y(i))\ge \delta n^{1/m}/k)\\&\qquad \le
 \sum_{n=1}^\infty \sum_{k=1}^\infty
 e^{-\lambda}\frac{\lambda^k}{k!} \sum_{i=1}^k \sup_{x\in X}P^x(\operatorname{dist}(Y(i-1),Y(i))\ge \delta n^{1/m}/k)\\
 &\qquad \le\sum_{n=1}^\infty \sum_{k=1}^\infty  e^{-\lambda}\frac{\lambda^k k}{k!\,\lambda}\sup_{x\in X}\int _{B(x,\delta n^{1/m}/k)^c}\varkappa(x,y)\,dy \\ &\qquad \le  \sum_{n=1}^\infty \sum_{k=1}^\infty  e^{-\lambda} \frac{\lambda^{k-1}}{(k-1)!}\, C\bigg(\frac k{\delta n^{1/m}}\bigg)^\alpha\\ &\qquad =\frac{e^{-\lambda} C}{\delta^\alpha}\bigg(\sum_{n=1}^\infty \frac1{n^{\alpha/m}}\bigg)\bigg(\sum_{k=1}^\infty\frac{\lambda^{k-1}k^\alpha}{(k-1)!}\bigg)<\infty,
    \end{align*}
  since $\alpha>m$.\qquad $\square$

  \section{Free Glauber dynamics as a scaling limit of free
  Kawasaki dynamics}\label{gggh}
Let  $\mu$ be a probability measure on $(\Gamma,{\cal B}(\Gamma))$. Assume
that,  for any $n\in\N$, there exists a non-negative measurable symmetric
function $k_\mu^{(n)}$ on $X^n$
 such
that, for any measurable symmetric function $f^{(n)}:X^n\to[0,\infty]$,
\begin{align} &\int_\Gamma \sum_{\{x_1,\dots,x_n\}\subset\gamma}
f^{(n)}(x_1,\dots,x_n)\,\mu(d\gamma)\notag\\&\qquad =\frac1{n!}\, \int_{X^n}
f^{(n)}(x_1,\dots,x_n) k_\mu^{(n)}(x_1,\dots,x_n)\,dx_1\dotsm
dx_n.\label{6t565r7}\end{align} Then, the functions $k_\mu^{(n)}$, $n\in\N$,
are called the correlation functions of the measure $\mu$.

Via a recursion formula, one can transform the correlation functions
$k_\mu^{\sk{n}}$ into the Ursell functions $u_\mu^{\sk{n}}$ and vice
versa, see e.g.\ \cite{Ru69}. Their relation is given by
\begin{equation}\label{urs_def}
k_\mu\sk{\eta}=\sum u_\mu\sk{\eta_1}\dotsm u_\mu\sk{\eta_j},\quad
\eta\in\Gamma_0,
\end{equation}
where
\[
\Gamma_0:=\skf{\eta\subset X: 1\le \md{\eta}<\infty },
\]
for any $\eta=\skf{x_1,\ldots,x_n}\in\Gamma_0$
\[
k_\mu\sk{\eta}:=k_\mu^{\sk{n}}\sk{x_1,\ldots,x_n}, \qquad
u_\mu\sk{\eta}:=u_\mu^{\sk{n}}\sk{x_1,\ldots,x_n} ,
\]
and the summation in (\ref{urs_def}) is over all partitions
of the set $\eta$ into nonempty mutually disjoint subsets
$\eta_1,\dots,\eta_j\subset\eta$ such that
$\eta_1\cup\dotsm\cup\eta_j=\eta$, $j\in\N$. Note that $k_\mu^{(1)}=u_\mu^{(1)}$.

Let now $X=\R^d$. We fix an arbitrary function $\xi\in S(\R^d)$ such that
 $\xi(-x)=\xi(x)$ for all $x\in\R^d$. Here, $S(\R^d)$
denotes the Schwartz  space of rapidly decreasing, infinitely differentiable
functions on $\R^d$. We define $$\varkappa(x,y):=\xi(x-y),\qquad x,y\in \R^d.$$
It can be easily checked that, by Theorem \ref{vgfuft}, the corresponding
Kawasaki dynamics exists as a Markov process
   ${\bf M}^{\mathrm K}$ on $(\Theta,{\cal B}(\Theta))$ with cadlag paths.

Let $\mu$ be a probability measure on $(\Theta,{\cal B}(\Theta))$ that
satisfies the following conditions:

(i) $\mu$ has correlation functions $(k_\mu^{(n)})_{n\in\N}$, and there exist
$0\le \gamma<1$ and $C>0$ such that
\begin{equation}\label{swaswea957}\forall n\in\N,\, \forall (x_1,\dots,x_n)\in (\R^d)^n:\quad
k_\mu^{(n)}(x_1,\dots,x_n)\le (n!)^\gamma C^n.\end{equation}

(ii) $\mu$ is translation invariant. In particular, the first correlation
function $k_\mu^{(1)}$ is a constant.

(iii) $\mu$ has decay of correlations in the following sense: for each $n\ge2$ and $1\le i\le n$
\begin{equation} \label{oiuyuy}u_\mu^{(n)}\left(\frac{x_1}\varepsilon,\dots,\frac{x_i}\varepsilon ,
x_{i+1},\dots,x_n\right)\to0\quad \text{as }\varepsilon\to0,
\end{equation}
 where the
convergence is in the $dx_1\dotsm dx_n$-measure on each compact set in $\R^d$.

For example, any double-potential Gibbs measure in the low activity-high
temperature regime satisfies the above assumptions, see \cite{Mi67,Ru69,Bro80}

Let us assume that the initial distribution of the Kawasaki dynamics is $\mu$.
We denote this stochastic process  by ${\bf M}_\mu^{\mathrm K}$. We scale this
dynamics as follows. Instead of the function $\xi$ used for the construction of
${\bf M}_\mu^{\mathrm K}$, use the function
$$\xi_\varepsilon(x):=\varepsilon^d\xi(\varepsilon x),\qquad x\in\R^d,$$
and denote the corresponding Kawasaki dynamics with initial distribution $\mu$
by ${\bf M}_{\mu,\,\eps}^{\mathrm K}$. We are interested in the limit of this
dynamics as $\eps\to0$.

As in subsec.~\ref{gyfhh}, we construct the Glauber dynamics ${\bf M}^{\mathrm
G}$ using $ a(x):=\la\xi\ra$ and $z=k^{(1)}_\mu$. Here and below, for any $f\in
L^1(\R^d,dx)$, we denote $$\la f\ra:=\int_{\R^d}f(x)\,dx.$$ We assume that the
initial distribution of the Glauber dynamics is $\mu$ and denote this
stochastic process  by ${\bf M}_{\mu}^{\mathrm G}$.

Below, we will use $\dd\Gamma$ to denote the space of multiple configurations
over $\R^d$ equipped with the vague topology, see e.g.\ \cite{Kal} for details.
We have $\Theta\in{\cal B}(\dd\Gamma)$.

\begin{th} Under the above assumptions, consider
${\bf M}_{\mu,\eps}^{\mathrm K}$, $\eps>0$, and ${\bf M}_{\mu}^{\mathrm G}$ as
stochastic processes taking values in $\dd\Gamma$. Then,  ${\bf
M}_{\mu,\,\eps}^{\mathrm K}\to {\bf M}_{\mu}^{\mathrm G}$ as $\eps\to0$ in the
sense of weak convergence of finite-dimensional distributions.
\end{th}

\noindent {\it Proof}. Let $(p_{t,\,\eps})_{t\ge0}$ denote the semigroup in
$L^2(\R^d,dx)$ with generator $$ (L_\eps
f)(x)=\int_{\R^d}(f(y)-f(x))\xi_\eps(x-y)\,dy.$$ Let $\mathcal F$ and ${\cal
F}^{-1}$ denote the Fourier transform and its inverse, respectively, which we
normalize so that they become unitary operators in $L^2(\R^d\to\C,dx)$. As
usual, we denote $\hat f:={\cal F}f$ and $\check f:={\cal F}^{-1}f$.

We easily have:\begin{equation}\label{tsrehy}
(p_{t,\,\eps}f)(x)=e^{-t\la\xi\ra}f(x)+ (K_{t,\,\eps}f)(x),
\end{equation}
where \begin{equation}\label{gutf}
(K_{t,\,\eps}f)(x):=\int_{\R^d}\eps^dG_{t}
(\eps(x-y))f(y)\,dy,\qquad t>0.\end{equation}
Here
\begin{equation}\label{jjijoih}
G_t(x):=e^{-t\la\xi\ra}(\exp[t(2\pi)^{d/2}\hat\xi]-1)\check{},\qquad
x\in\R^d.\end{equation} We note that, since $\xi\in S(\R^d)$, we have
$\hat\xi\in S(\R^d)$, and therefore $\exp[t(2\pi)^{d/2}\hat\xi]-1\in S(\R^d)$.
Hence, for each $t>0$, $G_t\in S(\R^d)$. Furthermore, since $\xi(-x)=\xi(x)$,
we get $G_t(x)=G_t(-x)$.

We fix any $n\in\N$, $0=t_0<t_1<t_2<\dots<t_n$, and
$\varphi_0,\varphi_1,\dots,\varphi_n\in C_0(\R^d)$ with $-1<\varphi_i\le0$,
$i=0,1\dots,n$. We denote by $({\bf P}_t^\eps(\gamma,\cdot))_{t\ge0,\,
\gamma\in\Theta}$ the transition semigroup of the $\eps$-Kawasaki dynamics for
$\eps>0$ and that of the Glauber dynamics for $\eps=0$. We also denote by
$(p^\eps_t(x,\cdot))_{t\ge0,\, x\in\R^d}$ the transition semigroup of the
one-particle $\eps$-dynamics, $\eps>0$.

By \eqref{6t565r7}, we have:
\begin{align}
&\int_\Theta \mu(d\gamma_0)\int_
\Theta {\bf P}^\eps_{t_1}(\gamma_0,d\gamma_1)\int_{\Theta}{\bf P}^\eps_{t_2-t_1}
(\gamma_1,d\gamma_2)\notag\\
&\qquad\times\dots\times\int_\Theta{\bf P}^\eps_{t_n-t_{n-1}}
(d\gamma_{n-1},d\gamma_n)\prod_{i=0}^n \exp[\la \log(1+\varphi_i),\gamma_i\ra]\notag\\
&\quad =\int_\Theta\mu(d\gamma)\prod_{x\in\gamma}
(1+\varphi_0(x))\int_{\R^d}p^\eps_{t_1}(x,dx_1)\int_{\R^d}p^\eps_{t_2-t_1}(x_1,dx_2)
\notag\\
&\qquad \times \dots\times \int_{\R^d}p^\eps_{t_n-t_{n-1}}(x_{n-1},dx_n)\prod_{i=1}^n (1+\varphi_i(x_i))\notag\\
&\quad =\int_\Theta\mu(d\gamma)\prod_{x\in\gamma} (1+g^\eps(x)),\notag\\
&=1+\sum_{m=1}^\infty \frac1{m!}\int_{(\R^d)^m}g^\eps(x_1)\dotsm g^\eps(x_m)
k_\mu^{(m)} (x_1,\dots,x_m)\,dx_1\dotsm dx_m.\label{kjghkug}
\end{align}
Here,
\begin{equation}\label{iguy} g^\eps(x):=\varphi_0(x)+f^\eps(x)+\varphi_0(x)f^\eps(x)\end{equation}
with \begin{align}f^\eps(x):&=\sum_{1\le i_1 <i_2<\dots<i_k\le n,\, k\ge1}
\int_{\R^d}p^\eps_{t_{i_1}}(x,dx_{1})\notag\\
&\quad \times\int_{\R^d}p^\eps_{t_{i_2}-t_{i_1}}
(x_{1},dx_{2})\times\dots\times p^\eps_{t_{i_k}-t_{i_{k-1}}}(x_{k-1},dx_k)
\varphi_{i_1}(x_1)\dotsm\varphi_{i_k}(x_k).\label{ytdstrswer5}\end{align}

We easily have from \eqref{iguy} and  \eqref{ytdstrswer5}: $$
\sup_{\eps>0}\int_{\R^d}|g^\eps(x)|\,dx<\infty.$$ Hence, by \eqref{swaswea957},
in oder to find the limit of \eqref{kjghkug} as $\eps\to0$, it suffices to find
the limit of each term in the sum.

By \eqref{tsrehy}--\eqref{jjijoih}, \eqref{iguy}, and \eqref{ytdstrswer5},
\begin{align}
g^\eps(x)&= \varphi_0(x)+(1+\varphi_0(x))\sum_{1\le i_1 <i_2<\dots<i_k\le n,\, k\ge1}(e^{-t_{i_1}\la\xi\ra}+K_{t_{i_1},\,\eps})M_{\varphi_{i_1}}\notag\\
&\quad\times
(e^{-(t_{i_2}-t_{i_1})\la\xi\ra}+K_{t_{i_2}-t_{i_1},\,\eps})M_{\varphi_{i_2}}\times\dots\times
(e^{-(t_{i_k}-t_{i_{k-1}})\la\xi\ra}+K_{t_{i_k}-t_{i_{k-1}},\,\eps})\varphi_{i_k},\label{iurf7r}
\end{align}
where $M_f$ denotes the operator of multiplication by a function $f$.

By \eqref{urs_def}, for each $m\in\N$,
\begin{align}&
\int_{(\R^d)^m}g^\eps(x_1)\dotsm g^\eps(x_m)
k_\mu^{(m)} (x_1,\dots,x_m)\,dx_1\dotsm dx_m\notag\\
&\qquad =\sum_{\{\eta_1,\dots,\eta_j\}}\prod_{i=1}^j
\int_{(\R^d)^{|\eta_i|}}g^{\eps}(x_1)\dotsm
g_\eps(x_{|\eta_i|})u_\mu^{(|\eta_i|)} (x_1,\dots,x_{|\eta_i|})\,dx_1\dotsm
dx_{|\eta_i|}, \label{kgyudc}
\end{align}
where the summation  is over all partitions $\{\eta_1,\dots,\eta_j\}$, $j\ge1$,
of the set $\{1,\dots,n\}$ into nonempty, mutually disjoint subsets
$\eta_1,\dots,\eta_j\subset\eta$.

We next have the following

\begin{lem}\label{uyrfydr}  Let the assumptions above be satisfied.
Let $k,n\in\N$, $k\le n$, $l_1,\dots,l_k\in\N$, $t^{(j)}_{i}>0$,
$j=1,\dots,l_i$, $i=1,\dots,k$. Let $f_i^{(1)},\dots,f_i^{(l_i)}\in C_0(\R^d)$,
$i=1,\dots,k$. Let $F:(\R^d)^n\to\R$ and $f_1,\dots,f_k:\R^d\to\R$ be
measurable and bounded, and $f_{k+1},\dots f_n\in C_0(\R^d)$. For any $\eps>0$,
set
\begin{align}
I_\eps:&=\int_{(\R^d)^n}dx_1\dotsm dx_n\, F(x_1,\dots,x_n) \prod_{i=1}^k \bigg(
f_i(x_i)\int_{\R^d} dx_i^{(1)}\eps^d G_{t^{(1)}_{i}}
(\eps(x_i-x_i^{(1)}))f_i^{(1)}(x_i^{(1)})\notag\\
&\quad\times
\int_{\R^d}dx_i^{(2)}\eps^dG_{t^{(2)}_{i}}(\eps(x_{i}^{(2)}-x_i^{(1)}))
f_i^{(2)}(x_i^{(2)})\notag\\ &\quad\times\dots\times
\int_{\R^d}dx_i^{(l_i)}\eps^d
G_{t^{(l_i)}_{i}}(\eps(x_i^{(l_i-1)}-x_i^{(l_i)}))f_i^{(l_i)}(x_i^{(l_i)})\bigg)\prod_{j=k+1}^n
f_j(x_j) .\label{iutguy}
\end{align}
We then have:

\rom{(i)} If  at lest one $l_i\ge2$, then $I_\eps\to0$ as $\eps\to0$.

\rom{(ii)} If $n\ge2$, $l_1=\dots=l_k=1$ and $F=u_\mu^{(n)}$, then
$I_\eps\to0$ as $\eps\to0$.

\rom{(iii)} If $l_1=\dots=l_k=1$, $F=1$, and at least one $f_i\in C_0(\R^d)$. Then. $I_\eps\to0$ as $\eps\to0$.

\rom{(iv)} If $l_1=\dots=l_k=1$, $f_1=\dots=f_k=1$, and $F=1$, then for each
$\eps>0$,$$ I_\eps=\left(\prod_{i=1}^k (1-e^{-t^{(1)}_{i}\la\xi\ra}) \la
f_i^{(1)}\ra\right)\times\prod_{j=k+1}^n\la f_j\ra.$$
\end{lem}

\noindent{\it Proof}. In the right hans side of \eqref{iutguy}, make the following change of variables
$$ x_i'=\eps(x_i),\ (x_i^{(j)})'=\eps x_i^{(j)},\quad j=1,\dots,l_i-1,\ i=1,\dots,n.$$
Then by the majorized convergence theorem and \eqref{oiuyuy}, we get the statements (i)--(iii). In the same way, in the case of (iv), we get
$$ I_\eps=  \left(\prod_{i=1}^k\la G_{t^{(1)}_i}\ra\la f_i^{(1)}\ra\right)\times\prod_{j=k+1}^n\la f_j\ra.$$
By  \eqref{tsrehy}, we have, for any $t>0$:
$$\la G_t\ra=1-e^{-t\la\xi\ra},$$ from where the lemma follows.\quad $\square$

By \eqref{urs_def},  \eqref{iurf7r}, \eqref{kgyudc}, and Lemma~\ref{uyrfydr},
and taking into account that each Ursell function $u_\mu^{(n)}$ is bounded, we
have:
\begin{align}& \int_{(\R^d)^m}g^\eps(x_1)\dotsm g^\eps(x_m)
k_\mu^{(m)} (x_1,\dots,x_m)\,dx_1\dotsm dx_m\notag\\&\quad \to
\sum_{l=0}^m\binom{m}{l } \int_{(\R^d)^{l}}
\prod_{j=1}^l\left(\varphi_0(x_j)+(1+\varphi_0(x_j)) \sum_{1\le i_1
<i_2<\dots<i_k\le n,\, k\ge1}
e^{-t_{i_k}\la\xi\ra}(\varphi_{i_1}\dotsm\varphi_{i_k})(x_j)
\right)\notag\\
&\quad\quad\times k_\mu^{(l)}(x_1,\dots,x_{l})\,dx_1\dotsm dx_{l} \notag\\
&\quad\quad\times\left(\sum_{1\le i_1 <i_2<\dots<i_k\le n,\, k\ge1}
k_\mu^{(1)}(1-e^{-t_{i_1}\la\xi\ra})e^{-(t_{i_k}-t_{i_1})\la\xi\ra}\la\varphi_{i_1}
\dotsm\varphi_{i_k}\ra\right)^{m-l}.
\end{align}
Therefore, by \eqref{6t565r7}, the left hand side of \eqref{kjghkug} converges
to
\begin{align}&\exp\left[\sum_{1\le i_1 <i_2<\dots<i_k\le n,\, k\ge1} k_\mu^{(1)}(1-e^{-t_{i_1}\la\xi\ra})e^{-(t_{i_k}-t_{i_1})\la\xi\ra}\la\varphi_{i_1}\dotsm\varphi_{i_k}\ra\right]\notag\\
&\quad\times\int_\Theta \exp[\la\log(1+\varphi_0),\gamma\ra]\prod_{x\in\gamma}
\left(1+\sum_{1\le i_1 <i_2<\dots<i_k\le n,\, k\ge1} e^{-t_{i_k}\la\xi\ra}
(\varphi_{i_1}\dotsm\varphi_{i_k})(x) \right)\,\mu(d\gamma)\label{raewaa}
\end{align}
as $\eps\to0$.

We next have the following

\begin{lem}\label{ertsresa} For each $\gamma\in\Theta$ and any $0<t_1<t_2,\dots<t_n$, $n\in\N$,
\begin{align}&\int_\Theta {\bf P}^0_{t_1}(\gamma,d\gamma_1)\int_{\Theta}{\bf P}^0_{t_2-t_1}
(\gamma_1,d\gamma_2)\times\dots\times\int_\Theta{\bf P}^0_{t_n-t_{n-1}}
(\gamma_{n-1},d\gamma_n)\prod_{i=1}^n \exp[\la \log(1+\varphi_i),\gamma_i\ra]\notag\\
&\quad =\exp\left[\sum_{1\le i_1 <i_2<\dots<i_k\le n,\, k\ge1}
k_\mu^{(1)}(1-e^{-t_{i_1}\la\xi\ra})e^{-(t_{i_k}-t_{i_1})\la\xi\ra}
\la\varphi_{i_1}\dotsm\varphi_{i_k}\ra\right]\notag\\
&\qquad\times \prod_{x\in\gamma}\left(1+\sum_{1\le i_1 <i_2<\dots<i_k\le n,\,
 k\ge1} e^{-t_{i_k}\la\xi\ra}(\varphi_{i_1}\dotsm\varphi_{i_k})(x)
\right).\label{uiiyutg}
\end{align}
\end{lem}

\noindent{\it Proof.} By Sec.~\ref{lkughihgh} and Subsec.~\ref{gyfhh}, for $n=1$, we have:
$$ \int_\Theta \exp[\la \log(1+\varphi_1),\gamma_1\ra]\,{\bf P}^0_{t_1}
(\gamma,d\gamma_1)=\exp\left[k_\mu^{(1)}(1-e^{-t_1\la\xi\ra})
\la\varphi_1\ra\right ]\prod_{x\in\gamma}(1+e^{-t_1\la\xi\ra}\varphi_1(x)),$$
which is \eqref{uiiyutg} in this case.

Now, assume that \eqref{uiiyutg}
 holds for $n\in\N$, and let us prove it for $n+1$. We then have:
 \begin{align}
 &\int_\Theta {\bf P}^0_{t_1}(\gamma,d\gamma_1)
 \int_{\Theta}{\bf P}^0_{t_2-t_1}
(\gamma_1,d\gamma_2)\times\dots\times\int_\Theta{\bf P}^0_{t_{n+1}-t_{n}}
(\gamma_{n},d\gamma_{n+1})\prod_{i=1}^{n+1} \exp[\la \log(1+\varphi_i),\gamma_i\ra]
\notag\\
&\quad=\int_\Theta {\bf
P}^0_{t_1}(\gamma,d\gamma_1)\exp[\log(1+\varphi_1),\gamma_1\ra]\notag\\
&\qquad\times \prod_{x\in\gamma_1}\left(1+\sum_{2\le i_1<\dots<i_k\le n+1,\,
k\ge 1}e^{-(t_{i_k}-t_1)\la\xi\ra }(\varphi_{i_1}\dotsm\varphi_{i_k})(x)\right)\notag\\
& \qquad \times \exp\left[\sum_{2\le i_1<\dots<i_k\le n+1,\, k\ge 1}
k_\mu^{(1)}(1-e^{-(t_{i_1}-t_1)\la\xi\ra })e^{-(t_{i_k}-t_{i_1})\la\xi\ra }\la
\varphi_{i_1}\dotsm
\varphi_{i_k}\ra \right]\notag\\
&=\quad \prod_{x\in\gamma} \left(e^{-t_1\la\xi\ra
}(1+\varphi_1(x))\left(1+\sum_{2\le i_1<\dots<i_k\le n+1,\, k\ge
1}e^{-(t_{i_k}-t_1)\la\xi\ra
}(\varphi_{i_1}\dotsm\varphi_{i_k})(x)\right)\right.\notag\\
&\qquad\left.\vphantom{\sum_{2\le i_1<\dots<i_k\le n+1,\, k\ge
1}e^{-(t_{i_k}-t_1)}} +1-e^{-t_1\la\xi\ra }
\right)\notag\\
&\qquad\times \int_\Theta\pi_{(1+e^{-t_1\la\xi\ra} )k_\mu^{(1)}}(d\gamma_1)
\prod_{x\in\gamma_1}(1+\varphi_1(x))\notag\\ &\qquad\times\left(1+\sum_{2\le
i_1<\dots<i_k\le n+1,\, k\ge
1}e^{-(t_{i_k}-t_1)\la\xi\ra }(\varphi_{i_1}\dotsm\varphi_{i_k})(x)\right)\notag\\
& \qquad \times \exp\left[\sum_{2\le i_1<\dots<i_k\le n+1,\, k\ge 1}
k_\mu^{(1)}(1-e^{-(t_{i_1}-t_1)\la\xi\ra })e^{-(t_{i_k}-t_{i_1})\la\xi\ra }\la
\varphi_{i_1}\dotsm
\varphi_{i_k}\ra \right]\notag\\
&\quad= \prod_{x\in\gamma}\left(1+\sum_{1\le i_1 <i_2<\dots<i_k\le n+1,\,
k\ge1} e^{-t_{i_k}\la\xi\ra}(\varphi_{i_1}\dotsm\varphi_{i_k})(x)
\right)\notag\\
&\qquad\times \exp\left[(1-e^{-t_1\la\xi\ra})k_\mu^{(1)}\vphantom{\sum_{2\le
i_1<\dots<i_k\le
n+1,\, k\ge 1}}\right.\notag\\
&\qquad\times \left.\left\la -1+(1+\varphi_1)\left(1+\sum_{2\le
i_1<\dots<i_k\le n+1,\, k\ge 1}e^{-(t_{i_k}-t_1)\la\xi\ra
}(\varphi_{i_1}\dotsm\varphi_{i_k})(x)\right) \right\ra \right.\notag\\
&\qquad\left.+\sum_{2\le i_1<\dots<i_k\le n+1,\, k\ge 1}
k_\mu^{(1)}(1-e^{-(t_{i_1}-t_1)\la\xi\ra })e^{-(t_{i_k}-t_{i_1})\la\xi\ra }\la
\varphi_{i_1}\dotsm \varphi_{i_k}\ra\right]\notag\\
&\quad= \prod_{x\in\gamma}\left(1+\sum_{1\le i_1 <i_2<\dots<i_k\le n+1,\,
k\ge1} e^{-t_{i_k}\la\xi\ra}(\varphi_{i_1}\dotsm\varphi_{i_k})(x)
\right)\notag\\
&\qquad \times\exp\left[\sum_{1\le i_1 <i_2<\dots<i_k\le n+1,\, k\ge1}
k_\mu^{(1)}(1-e^{-t_{i_1}\la\xi\ra})e^{-(t_{i_k}-t_{i_1})\la\xi\ra}
\la\varphi_{i_1}\dotsm\varphi_{i_k}\ra\right].\notag
\end{align}
Thus, by induction, the lemma is proved. \quad $\square$

Let $0\le t_0<t_1<\dots<t_n$, and let $\mu^\eps_{t_0,t_1,\dots,t_n}$,
$\eps\ge0$, denote the joint distribution of the  process ${\bf
M}_{\mu,\,\eps}^{\mathrm K}$  at times $t_0,t_1,\dots,t_n$ for $\eps>0$, and
respectively that of the process ${\bf M}_{\mu}^{\mathrm G}$ for $\eps=0$.

 By
\eqref{raewaa} and Lemma \ref{ertsresa}, for any
$\varphi_0,\varphi_1\dots,\varphi_n\in C_0(\R^d)$, $\varphi_i\le0$,
$i=0,1,\dots,n$,
\begin{align}
&\int_{\dd\Gamma{}^{n+1}} \prod_{i=0}^n \exp[\la \varphi_i,\gamma_i\ra]\,
d\mu^\eps_{t_0,t_1,\dots,t_n}(\gamma_0,\gamma_1,\dots,\gamma_n)
\notag\\
&\quad\to \int_{\dd\Gamma{}^{n+1}} \prod_{i=0}^n \exp[\la
\varphi_i,\gamma_i\ra]\,
d\mu^\eps_{t_0,t_1,\dots,t_n}(\gamma_0,\gamma_1,\dots,\gamma_n)
 \label{erweawea}
\end{align}
as $\eps\to0$. By \eqref{erweawea}, we, in particular, get, for any $t>0$ and
any $\varphi\in C_0(\R^d)$, $\varphi\le0$,
$$ \int_{\dd\Gamma}\exp[\la
\varphi,\gamma\ra]\,d\mu_t^\eps(\gamma)\to \int_{\dd\Gamma}\exp[\la
\varphi,\gamma\ra]\,d\mu_t^0(\gamma)$$ as $\eps\to0$. Hence, by \cite[Theorem~4.2]{Kal},
$\mu_t^\eps\to\mu_t^0$ weakly in ${\mathcal M}(\dd\Gamma)$ as $\eps\to0$. Here,
${\mathcal M}(\dd\Gamma)$ denotes the space of probability measures on
$\dd\Gamma$, see e.g. \cite{Par} for details. Therefore, the set
$\{\mu_t^\eps\mid 0<\eps\le1\}$ is tight in ${\mathcal M}(\dd\Gamma)$ . This
implies that, for any $0\le t_0<t_1<\dots<t_n$, the set
$\{\mu_{t_0,t_1,\dots,t_n}^\eps\mid 0<\eps\le1\}$ is tight in ${\mathcal
M}(\dd\Gamma{}^{n+1})$. Hence, by \eqref{erweawea},
$\mu_{t_0,t_1,\dots,t_n}^\eps\to\mu_{t_0,t_1,\dots,t_n}^0$ weakly in ${\mathcal
M}(\dd\Gamma{}^{n+1})$ as $\eps\to0$. Thus, the the theorem is proved. \quad
$\square$

  \begin{center}
{\bf Acknowledgements}\end{center}

 The authors acknowledge the financial support of the SFB 701 `` Spectral
structures and topological methods in mathematics''.

\end{document}